 \documentclass[10pt]{article}
\font\smallit=cmti10
\font\smalltt=cmtt10

\usepackage{amssymb,latexsym,amsmath,epsfig,amsthm} 

\makeatletter

\renewcommand\section{\@startsection {section}{1}{\z@}
{-30pt \@plus -1ex \@minus -.2ex}
{2.3ex \@plus.2ex}
{\normalfont\normalsize\bfseries\boldmath}}

\renewcommand\subsection{\@startsection{subsection}{2}{\z@}
{-3.25ex\@plus -1ex \@minus -.2ex}
{1.5ex \@plus .2ex}
{\normalfont\normalsize\bfseries\boldmath}}

\renewcommand{\@seccntformat}[1]{\csname the#1\endcsname. }

\makeatother

\newtheorem{theorem}{Theorem}
\newtheorem{lemma}{Lemma}
\newtheorem{proposition}{Proposition}
\newtheorem{corollary}{Corollary}
\newtheorem{fact}{Fact}

\theoremstyle{definition}

\newtheorem{remark}{Remark}

\usepackage{amsmath}
\usepackage{amsthm}
\usepackage{amsfonts}
\usepackage{amssymb}
\usepackage[all]{xy}
\usepackage{alltt}
\usepackage{enumitem}
\usepackage{stmaryrd}
\usepackage{comment}
\usepackage{xspace}
\usepackage{comment}
\usepackage{diagbox}
\usepackage{adjustbox}
\usepackage[thinlines]{easytable}
\usepackage{verbatim}
\usepackage{MnSymbol} 

\newcommand{\R}{\ensuremath{\mathbb{R}}}

\newcommand{\rar}{\rightarrow}
\newcommand{\intsn}{\mathbb{Z}_n}
\newcommand{\ints}{\mathbb{Z}}

\def \n{\noindent }
\def \bs{\bigskip}
\def \R{\mathbb R}

\def \G{\mathcal G}
\def \R{\mathcal R}
\def \I{\mathcal I}

\def \al{\alpha}

\def \x{{\bf x}}
\def \y{{\bf y}}

\def \e{{\bf e}}
\def \cv{{\bf c}}

\DeclareMathOperator{\rank}{rank}

\usepackage[
 colorlinks=true,urlcolor=blue,linkcolor=blue,citecolor=blue
]{hyperref}

\begin{document}  

\begin{center}
\uppercase{\bf Zero-sum-free tuples and hyperplane arrangements}
\vskip 20pt
{\bf Sunil K. Chebolu\footnote{The first author is supported by Simons Foundation: Collaboration Grant for Mathematicians (516354).}}\\
{\smallit Department of Mathematics, Illinois State University, Normal, IL, USA}\\
{\tt schebol@ilstu.edu}\\  
\vskip 10pt
{\bf Papa A. Sissokho}\\
{\smallit Department of Mathematics, Illinois State University, Normal, IL, USA}\\
{\tt psissok@ilstu.edu}\\ 
\end{center}
\vskip 20pt

\vskip 30pt


 %
\centerline{\bf Abstract}

\noindent
A vector $(v_{1}, v_{2}, \cdots, v_{d})$ in $\intsn^{d}$ is 
said to be a zero-sum-free $d$-tuple if there is no non-empty subset of its components whose sum is zero in $\mathbb{Z}_n$. We denote the cardinality of this collection by $\alpha_n^d$. We let $\beta_n^d$ denote the cardinality of the set of zero-sum-free tuples in $\intsn^{d}$ where $\gcd(v_1, \cdots,v_d, n) = 1$. We show that $\alpha_n^d=\phi(n)\binom{n-1}{d}$ when $d > n/2$, and in the general case, we prove recursive formulas, divisibility results, bounds, and asymptotic results for $\alpha_n^d$ and $\beta_n^d$. In particular, $\alpha_n^{n-1} = \beta_n^1= \phi(n)$, suggesting that these sequences can be viewed as generalizations of Euler's totient function. 
We also relate the problem of computing $\alpha_n^d$ to counting points in the complement of a certain hyperplane arrangement defined over $\mathbb{Z}_n$. It is shown that the hyperplane arrangement's characteristic polynomial captures $\alpha_n^d$ for all integers $n$ that are relatively prime to some determinants. We study the row and column patterns in the numbers $\alpha_n^{d}$. We show that for any fixed $d$, $\{\alpha_n^d \}$ is asymptotically equivalent to $\{ n^d\}$. We also show a connection between the asymptotic growth of $\beta_n^d$ and the value of the Riemann zeta function $\zeta(d)$. Finally, we show that $\alpha_n^d$ arises naturally in the study of Mathieu-Zhao subspaces in products of finite fields.

\pagestyle{myheadings}
\markright{\smalltt INTEGERS: 21 (2021)\hfill}
\thispagestyle{empty}
\baselineskip=12.875pt
\vskip 30pt 

\section{Introduction}\label{sec:introduction}

Let $\mathbb{Z}_n$ denote the ring of integers modulo $n$. A vector $(v_{1}, v_{2}, \cdots, v_{d})$ in $\intsn^{d}$ is 
said to be a zero-sum-free $d$-tuple if there is no non-empty subset of its components whose sum is zero in $\mathbb{Z}_n$. We denote this collection by $\G_n^d$, and our main goal is to compute $\alpha_n^d$, the cardinality of $\G_n^d$.

Zero-sum-free sequences over an (additively written) abelian group $H$ have been extensively studied (e.g., see~\cite{GGp, SC, VP}).
In particular, this study is usually done for the concept of {\em minimal} zero-sum sequence over 
$H$, which is a sequence $S$ (where order does not matter) of elements of $H$ that sum to the identity element $0$ in $H$ and such that no proper subsequence of $S$ adds to $0$ (see the survey papers~\cite{C, GGs}). A central question in this area is the determination of the {\em Davenport constant of $H$}, which is the smallest positive integer $t$ such that any sequence of length $t$ contains a non-empty subsequence with sum $0$. 
Although in this paper, we work with $d$-tuples instead of sequences, the two are related, as the reader will see in what follows. In particular, we draw upon the work on minimal zero-sum sequences. However, our initial motivation for studying zero-sum-free $d$-tuples came from a problem in algebra involving Mathieu-Zhao subspaces, which we explain in the last section.

In this paper, we provide some concrete approaches to the problem of computing $\alpha_n^d$. Closely related to this is $\beta_n^d$, the cardinality of the set of irreducible zero-sum-free $d$-tuples in $\mathbb{Z}_n^d$, where $(x_1, \cdots, x_d)$ is irreducible if $\gcd(x_1, \cdots,x_d, n) = 1$. As we will see, computing one of these is equivalent to computing the other.
Our first computational technique is based on group actions, and the second one uses hyperplane arrangements.  
Before explaining our main results, we introduce some notation.
For any prime $p$, let $\mathbb{F}_p$ denote the field with $p$ elements, and let $\mathbb{F}_p^d$ denote the vector space of dimension $d$ over $\mathbb{F}_p$.  
 For any positive integer $n$, let $\phi(n)$ denote the Euler totient function, which is the number of positive integers $t$ such that $t\leq n$ and $\gcd(t,n)=1$. The set $\G_n^{d}$ carries a natural action of $Aut(\mathbb{Z}_n)$, where each automorphism of $\mathbb{Z}_n$ acts component-wise. Using this action and the characterization of minimal zero-sum sequences, we investigated the numbers $\alpha_n^d$ and $\beta_n^d$. We restrict to $1 \le d < n$ because it is easy to see that $\alpha_n^d = 0$ for all $d \ge n$; see Fact \ref{fact:1}. Our first main theorem is: 

\begin{theorem}\label{thm:1.1}
Let $n\geq 3$ and $1\leq d<n$ and consider the natural action of $Aut(\intsn)$ on $\G_n^d$. Then we have the following.
\begin{enumerate}
\item If $d > n/2$, then $\alpha_n^d = \phi(n) {n-1 \choose d }$.
\item If $d\geq \max\{m:\, 1\leq m<n,\, m\mid n\}$, then $\phi(n)$ divides $\alpha_n^d$.
 \item $\al_n^d=\sum\limits_{m\mid n,\, m \geq 1}\beta_m^d\, \mbox{ and }\, \beta_n^d=\sum\limits_{m\mid n,\,m\geq 1}\mu\left(\frac{n}{m}\right)\al_{m}^d$,
 where $\mu$ is the M\"obius function.
\item $\phi(n)$ divides $\beta_n^d$.
\end{enumerate}
\end{theorem}
\begin{remark}
We prove parts A, B, C,  and D of Theorem~\ref{thm:1.1} in Theorems~\ref{savchev-chen}, \ref{thm:cha}(d), 
\ref{thm:cha}(c), and \ref{thm:cha}(b), respectively.
\end{remark}

Setting $d= n-1$ in part A, we see that $\alpha_n^{n-1} = \phi(n)$. Thus, we can think of $\alpha_n^d$ as a generalization of Euler's totient function $\phi(n)$ to two variables. This suggests that computing $\alpha_n^d$ can be as hard as computing $\phi(n)$. In fact, $\alpha_n^d$ fails to have the nice properties that $\phi(n)$ has. For instance, while $\phi(n)$ is a multiplicative function, $\alpha_n^d$ is not; see Corollary \ref{cor:rec}, part $(iii)$.
Note that part C says that determining $\{\al_n^d\}$ is equivalent to determining $\{\beta_n^d\}$. We also obtain recursive formulas and bounds for $\alpha_n^d$ and $\beta_n^d$; see Corollary~\ref{cor:rec} and Proposition \ref{bounds}.

Our next approach is via hyperplane arrangements. Note that $\G_n^d$ is the complement, in $\mathbb{Z}_n^d$, of the union of the $2^d -1$ hyperplanes defined by the equations $\sum_{i \in S} x_i = 0$, where $S$ ranges over the set of all non-empty subsets of $\{ 1, 2, \ldots, n\}$. So we can compute its cardinality using the exclusion-inclusion principle and techniques from hyperplane arrangements. For instance, using exclusion-inclusion, the first non-trivial formula is:
\[ \alpha_n^3 = n^3 - 7n^2 + 15n-10+ \frac{1+(-1)^{n-1}}{2} \text{ for all } n \ge 3.\]
Using this approach, we rediscover the characteristic polynomial associated with Hyperplane arrangements, which is known to capture $\alpha_p^d$ when $p$ is a sufficiently large prime number. We generalize this result and get information on the coefficients of these polynomials. $H_d$ will denote a $d \times (2^d-1)$ matrix whose columns are all the non-zero binary vectors in $\mathbb{F}_p^d$.
\begin{theorem}\label{thm:1.2} 
Let $d$ be a positive integer. 
\begin{enumerate} 
\item
There exists a monic polynomial $f_d(x)$ of degree $d$ with integer coefficients such that $\alpha_n^d = f_d(n)$ for all $n$ that are relatively prime to the determinant of any $d \times d$ binary matrix. 
\item If $\text{gcd}\left(n, \lceil d^{d/2}\rceil!\right) = 1$, then $\alpha_n^d = f_d(n)$.
\item The coefficient of $x^i$ ($0 \le i \le d-1$) in $f_d(x)$ is given by $\sum_{j=1}^{2^d-1} \;(-1)^j m(j, i)$, where 
$m(j, i)$ is the number of subsets of $j$ columns of $H_d$ that span a $d-i$ dimensional subspace in
 $\mathbb{F}_p^d$.
\end{enumerate}
\end{theorem}
\begin{remark}
We prove parts A, B, and C of  Theorem~\ref{thm:1.2} in Theorem~\ref{thm:prime-ext}, Corollary~\ref{gcdlemma}, and  Proposition~\ref{prop:coeff}, respectively.
\end{remark}

 We computed values of $\alpha_n^d$ using SageMath for various values of $n$ and $d$ and organized these values into a table where the rows correspond to values when $n$ is fixed and the columns to values when $d$ is fixed. Based on this data, we formulated row and column hypotheses (see Section~\ref{sec:conj}). Row hypothesis states that for all $n$, $\alpha_n^i$ is an increasing sequence for $0 < i < n/2$ and the column hypothesis states that for each $d$, $\alpha_i^d$ is an increasing sequence. We found counterexamples to both these hypotheses. This gave another hypothesis called the eventual column hypothesis, which states that for each $d$, $\alpha_n^d$ is increasing sequence for $n$ large enough. These hypotheses led us to study the asymptotic growth of $\{\alpha_n^d\}$ and $\{ \beta_n^d\}$. Our work on these hypotheses and asymptotic results can be summarized in the following theorem. 
\begin{theorem}\label{thm:1.3}
 Let $n\geq 3$ and $1< d<n$. 
\begin{enumerate}
\item There are infinitely many positive integers $d$ such that $\alpha_{n+1}^d < \alpha_n^d$ for some $n$ that depends on $d$. 
\item For every positive integer $d$, the sequence $\{ \alpha_n^d\}$ is an increasing sequence if $n$ is large enough.
\item For any fixed $d$, $\alpha_n^d/n^d \rightarrow 1$ as $n \rightarrow \infty$.
\item We have $1=\limsup_n \frac{\beta_n^d}{n^d} \ge \liminf_n \frac{\beta_n^d}{n^d} \; \ge \; \frac{1}{d!\; \zeta(d)} \ge 0 $, where $\zeta(s)$ is the Riemann zeta function. Moreover, when $d \ge 2$, then $\zeta(d) > 0$, and the sequence $\{ \beta_n^d \}$ is asymptotically bounded above and below by the sequence $\{ n^d \}$.
\end{enumerate}
\end{theorem}
\begin{remark}
We prove parts A, B, C of Theorem~\ref{thm:1.3} in Proposition~\ref{prop:nocol}, Theorem~\ref{thm:eventual-col}, and  Theorem~\ref{thm:asymptotic-m}(ii), respectively. Part D is proved by combining Theorem~\ref{thm:limsup}, Proposition~\ref{prop:zeta}, and Theorem~\ref{thm:bounded}.
\end{remark}

Part A tells us that while there are infinitely many counterexamples to the column hypothesis, part $B$ means that 
the eventual column hypothesis is true. Part C tells us that for any fixed $d$, the sequence $\{ \alpha_n^d \}$ is asymptotically equivalent to the sequence $\{ n^d\}$.
Note that when $d =1$, $\beta_n^1 = \phi(n)$, and item D can be viewed as an extension of the known properties (see \cite{TA})  of $\phi(n)$:
\[1 = \limsup_n \phi(n)/n > \liminf_n \phi(n)/n = 0.\]
Thus, $\beta_n^d$ can be viewed as yet another generalization for $\phi(n)$.

One can also approach counting points in the complement of a hyperplane arrangement using cohomological methods; see \cite{AT}.

\section{Counting zero-sum-free $d$-tuples} \label{counting}
Our goal is to compute $\alpha_n^d$, the cardinality of the collection $\G_n^d$ of zero-sum-free $d$-tuples in $\intsn^{d}$. We begin with an observation.

\begin{fact}\label{fact:1}
Let $d$ and $n$ are positive integers. Then $\al_n^d=0$ if and only if $d \ge n$.
\end{fact}
  
This follows from the well-known fact that $a_1,\ldots,a_n$ of (not necessarily distinct) elements of $ \mathbb{Z}_n$ contains a subsequence whose sum is $0\in \mathbb{Z}_n$, and that $(1, 1, \cdots, 1)$ is a zero-sum-free $d$-tuple in $\mathbb{Z}_n^{d}$ whenever $d< n$.
 So throughout we will assume that $d < n$. Moreover, since $\alpha_2^d = 1$ when $d=1$, and $0$ otherwise, we may further assume $n \ge 3$. 
\subsection{Action of $Aut(\intsn)$ on $\G_n^d$} \label{action}
Let $Aut(\ints_{n})$ denote the automorphism group of $\ints_{n}$. Then 
$Aut(\ints_{n})$ acts naturally on $\G_n^{d}$ component-wise. In this section and beyond, we will
 identify $\ints_{n}$ with the set $\{ 0, 1, \cdots, n-1 \}$, and $Aut(\ints_{n})$ with the set $\{k:\; 1\leq k\leq n-1\mbox{ and } \gcd(k,n)=1\}$.

A zero-sum free $d$-tuple $\x=(x_1,\ldots,x_d)\in \G_n^d$ is called {\em irreducible} if $\gcd(\x,n):=\gcd(x_1,\ldots,x_d,n)=1$.
Otherwise, $x$ is called {\em reducible}. Let $\R_n^d$ (resp. $\I_n^d$) denote the sets of reducible (resp. 
irreducible) zero-sum free $d$-tuple in $\G_n^d$. Then, 
\[\R_n^d\cap\I_n^d=\emptyset\mbox{ and } \G_n^d=\R_n^d\cup\I_n^d.\] 

 The following theorem shows that computing the sequence $\{\al_n^d\}$ of the number of zero-sum-free 
$d$-tuples is equivalent to computing the sequence $\{\beta_n^d\}$ of the number of {\em irreducible} zero-sum-free $d$-tuples. 

\begin{theorem}\label{thm:cha}
Let $n\geq 3$ and $1\leq d<n$ and consider the action of $Aut(\intsn)$ on $\G_n^d$.

\n $(a)$ For any $k\in Aut(\intsn)$ and any $\x\in\G_n^d$, we have $k\in Stab(\x)$ if and only if $\frac{n}{\gcd(k-1,n)}$ divides $\gcd(\x,n)$. In particular, if $\gcd(\x,n)=1$, then $Stab(\x)=\{1\}$ for any $\x\in \G_n^d$.

\n $(b)$ The number of orbits of the restricted action of $Aut(\intsn)$ on $\I_n^d$ is $\beta_n^d/\phi(n)$. Thus, $\phi(n)$ divides $\beta_n^d$.
 
\n $(c)$ Determining $\{\al_n^d\}$ is equivalent to determining $\{\beta_n^d\}$. More precisely, we have 
\[\al_n^d=\sum\limits_{m\mid n,\, m \geq 1}\beta_m^d\, \mbox{ and }\, \beta_n^d=\sum\limits_{m\mid n,\,m\geq 1}\mu\left(\frac{n}{m}\right)\al_{m}^d,\]
 where $\mu$ is the M\"obius function.
 
\n $(d)$ If $d\geq \max\{m:\; \mbox{$1\leq m<n$ and $m\mid n$}\}$, then  
$\alpha_n^d=\beta_n^d$ and $\phi(n)$ divides $\alpha_n^d$. Moreover, the conclusion 
of this statement holds if $d\geq n/2$.
 
\end{theorem}
\begin{proof} 
To prove  $(a)$, let $k\in Aut(\intsn)$ and $\x\in\G_n^d$. Then 
$k\in {\rm Stab}(\x)$ is equivalent to $kx_i\equiv x_i$ for $1\leq i\leq d$, which 
in turn is equivalent to  $n$ divides $(k-1)x_i$ for $1\leq i\leq d$. The latter statement is 
equivalent to $\frac{n}{\gcd(k-1,n)}$ divides $x_i$ for $1\leq i\leq d$, which is true if and only if 
\[\frac{n}{\gcd(k-1,n)} \mbox{ divides } \gcd(\x,n),\]
where $\gcd(\x,n)=\gcd(x_1,\ldots,x_d,n)$. This proves the ``if and only if'' statement in part $(a)$.
Finally, if $\gcd(\x,n)=1$, then $\frac{n}{\gcd(k-1,n)}$ divides $\gcd(x_1,\ldots,x_d,n)$, which is equivalent to  
$\frac{n}{\gcd(k-1,n)}$ divides $1$. The latter statement is equivalent to $\gcd(k-1,n)=n$, which holds if and only 
if $k=1$.

\bs\n To prove $(b)$, note that if $\gcd(x_1,\ldots,x_d,n) = 1$, then $\gcd(k x_1,\ldots, k x_d,n) = 1$ for any integer $k$ such that  $1 \le k \le n-1$ with $\gcd(k, n) = 1$. So, we can consider the restricted action of $Aut(\intsn)$ on $\I_n^d$. If $\x\in\I_n^d$, then by definition, we have $\gcd(\x,n)=1$. Then it follows from part $(a)$ that $Stab(\x)=\{1\}$. Thus, the size of the orbit of $\x\in \I_n^d$ is equal to $|Aut(\intsn)|/|{\rm Stab}(\x)|=\phi(n)/1=\phi(n)$. Then, the number of orbits of the restricted action of $Aut(\intsn)$ on $\I_n^d$ is $\beta_n^d/\phi(n)$. This proves $(b)$.

\bs\n For the proof of $(c)$, note that since $\al_{n}^d=|\G_n^d|=|\R_n^d|+\beta_n^d$, it suffices to show that 
\[|\R_n^d|=\sum_{m\mid n\atop 1<m<n} \beta_m^d.\]

Consider the multiset $\I=\bigsqcup\limits_{m\mid n\atop 1<m<n} \I_m^d$, where we use the symbol $\sqcup$ to denote the ``disjoint union''. 
Define the function $f:\; \R_{n}^d\to \I$ by  
\[f(\x)=\frac{1}{\gcd(\x,n)}\,\x=\left(\frac{x_1}{\gcd(\x,n)},\ldots,\frac{x_d}{\gcd(\x,n)}\right) \mbox{ and } f(\x)\in \I_m^d,\mbox{ where $m=\frac{n}{\gcd(\x,n)}$}.\]
In particular, $f(\x)$ and $f(\x')$ may be equal in value but still be considered different in $\I$ if they belong to $\I_m^d$ and $\I_{m'}^d$, respectively, such that 
\[m=\frac{n}{\gcd(\x,n)}\not= \frac{n}{\gcd(\x',n)}=m'.\]
We claim that $f$ is a bijection.
First, we prove that $f$ is onto. If $\y=(y_1,\ldots,y_d)\in \I$, then there exits $m$, with $1<m<n$, such that $m\mid n$ and 
$\y\in \I_{m}^d$. Let $\x=\frac{n}{m}\y$. Since $\y$ is irreducible, it follows that $\gcd(\y,n)=1$, which implies that 
\[\gcd(\x,n)=\gcd\left(\frac{n}{m}\y,n\right)=\frac{n}{m}.\]
Thus,
\[f(\x)=f\left(\frac{n}{m}\,\y\right)=\frac{1}{\gcd(\frac{n}{m}\,\y,n)}\left(\frac{n}{m}\,\y\right)=
 \frac{m}{n}\,\left(\frac{n}{m}\,\y\right)=\y,\] 
 showing that $f$ is onto.
To show that $f$ is one-to-one, let $\x,\,\x'\in \R_n^d$ be such that $f(\x)=f(\x')$. Since $\x$ and $\x'$ are reducible, it follows that $\gcd(\x,n)>1$ and $\gcd(\x',n)>1$. Thus, 
\[ m=\frac{n}{\gcd(\x,n)}<n \mbox{ and } m'=\frac{n}{\gcd(\x',n)}<n.\]
Moreover, since $\x$ and $\x'$ are zero-sum free, their entries are less than $n$. Thus, $\gcd(\x,n)<n$ and $\gcd(\x',n)<n$, which imply that  
\[ 1<m<n \mbox{ and } 1<m'<n.\]
If $m\not=m'$, then $f(\x)=\frac{1}{m}\x\in \I_m^d$ and $f(\x')=\frac{1}{m}\x'\in \I_{m'}^d$ are different copies of the same element in $\I$, which would contradict the hypothesis $f(\x)=f(\x')$. Thus, $m=m'$, which yields $f(\x)=f(\x')$.
Then, 
\[\frac{1}{\gcd(\x,n)}\,\x=\frac{1}{\gcd(\x',n)}\,\x',\]
which implies that 
\[\frac{m}{n}\,\x=\frac{m'}{n}\,\x', \mbox{ or, equivalently,  }\x=\x'.\]
This completes the proof that $f$ is a bijection.  Thus, 
\begin{equation}\label{eq:rcs}
\al_n^d=|\R_{n}^d|+ \beta_n^d=\sum_{m\mid n\atop 1<m<n} \beta_m^d+ \beta_n^d=\sum_{m\mid n,\,m>1} \beta_m^d=\sum_{m\mid n,\, m\geq 1} \beta_m^d,
\end{equation}
where the last equality holds since $\beta_1^d\leq \al_1^d$, and $\al_1^d=0$ by Fact~\ref{fact:1}.  
The resulting formula for $\beta_n^d$ follows directly from the above formula for $\al_n^d$ and 
M\"obius inversion formula. This completes the proof of $(c)$.

\bs\n Finally, we  prove $(d)$. If $d\geq \max\{m:\; \mbox{$1\leq m<n$ and $m\mid n$}\}$, then we know from Fact~\ref{fact:1} that $\al_{m}^d=0$ (and, thus, $\beta_m^d=0$) for all divisors $m$ of $d$ with $1\leq m<n$. Thus, it follows 
from statements $(c)$ and $(a)$ that $\alpha_n^d=\beta_n^d$ and $\phi(n)$ divides $\alpha_n^d$. For the second 
part of statement $(d)$, note that if $m\mid n$ and $m<n$, then $m\leq n/2$. 
\end{proof}
%


The following result shows more specific recursive formulas for $\al_n^d$ and $\beta_n^d$.
\begin{corollary}\label{cor:rec}
Let $n$ be an integer such that $n\geq 3$ and let $p_1,\ldots,p_t$ be the distinct prime divisors of $n$.
For any subset $S\subseteq \{p_1,\ldots,p_t\}$, define $n_S:=\frac{n}{\prod_{p\in S}p}$, where $\prod_{p\in S}p=1$ if $S=\emptyset$. Further assume that $1\leq d<n$.

\n $(i)$ Then, 

\[\al_n^d= \beta_n^d+\sum_{\emptyset\neq \,S\,\subseteq \{p_1,\ldots,p_t\}}(-1)^{|S|+1} \al_{n_S}^d\,
\mbox{ and } \, \beta_n^d=\sum_{S\subseteq \{p_1,\ldots,p_t\}}(-1)^{|S|} \al_{n_S}^d.\]
In particular,
\[ \phi(n)=\beta_n^1=\sum_{S\subseteq \{p_1,\ldots,p_t\}}(-1)^{|S|} (n_S-1).\]

\n $(ii)$ If $n=p^t$ for a prime $p$ and a positive integer $t$, then $\al_{p^{t}}^d=\beta_{p^t}^d+\al_{p^{t-1}}^d$.

\n $(iii)$ If $n=p_1p_2$ for two distinct primes $p_1$ and $p_2$, then $\al_{p_1p_2}^d=\beta_{p_1p_2}^d+\al_{p_1}^d+\al_{p_2}^d$.
\end{corollary}
\begin{proof}\
To prove $(i)$, note that since $\{p_1,\ldots,p_t\}$ is the set of distinct prime divisors of $n$, there exist positive integers $r_i$
such that $n=\prod_{i=1}^t p_i^{r_i}$ is a prime factorization of $n$. For $1\leq i\leq t$, we define 
\[n_i=\frac{n}{p_i} \mbox{ and } A_i=\{m:\; m\mid n_i,\, m\geq 1\}.\]
Then for any divisor $m$ of $n$ with $1\leq m<n$, there exits some $i$ such that $m$ divides $n_i$.
Thus, 
\begin{equation}\label{eq:inex}
\{m:\, m\mid n,\, 1\leq m<n\}=\bigcup_{i=1}^t \{m:\; m\mid n_i,\, m\geq 1\}=A_1\cup \ldots\cup A_t.
\end{equation}
Observe that for any subset $\{A_{i_1},\ldots, A_{i_k}\}\subseteq \{A_1,\ldots,A_t\}$, we have 
\begin{align}\label{eq:int}
A_{i_1}\cap \ldots\cap A_{i_k}
&=\left\{m:\; m\mid n_{i_j}\mbox{ for all $1\leq j\leq k$},\, m\geq 1\right\} \cr
&=\left\{m:\; m\mid \gcd(n_{i_1},\ldots,n_{i_k}),\, m\geq 1\right\}\cr
&=\left\{m:\; m\mid \gcd\left(\frac{n}{p_{i_1}},\ldots, \frac{n}{p_{i_k}}\right),\, m\geq 1\right\}\cr
&=\left\{m:\; m\,\big\vert\, \frac{n}{p_{i_1}\ldots p_{i_k}},\, m\geq 1\right\}\cr
&=\left\{m:\; m\mid n_S,\, S=\{p_{i_1},\ldots,p_{i_k}\},\, m\geq 1\right\},
\end{align}
where $n_S=\frac{n}{\prod_{p\in S}p}= \frac{n}{p_{i_1}\ldots p_{i_k}}$. 
Then by combining Theorem~\ref{thm:cha}$(c)$ with the relations in~\eqref{eq:inex} and~\eqref{eq:int}, and the inclusion-exclusion applied to $A_1\cup \ldots\cup A_t$, we obtain
\begin{align}
\al_n^d=\sum_{m\mid n,\, n\geq 1} \beta_m^d
&=\beta_n^d+\sum_{m\mid n,\,1\leq m<n} \beta_m^d\cr
&=\beta_n^d+\sum_{m\in A_1\cup \ldots\cup A_t} \beta_m^d\cr
&=\beta_n^d+\sum_{k=1}^t\;\sum_{1\leq i_1<\ldots<i_k\leq t} (-1)^{k+1} \sum_{m\in A_{i_1}\cap \ldots\cap A_{i_k}} \beta_m^d\cr
&=\beta_n^d+\sum_{\emptyset\neq \,S\,\subseteq \{p_1,\ldots,p_t\}} (-1)^{|S|+1}\sum_{m\mid n_S,\, n\geq 1} \beta_m^d\cr
&=\beta_n^d+\sum_{\emptyset\neq \,S\,\subseteq \{p_1,\ldots,p_t\}} (-1)^{|S|+1} \al_{n_S}^d.
\end{align}
The formula for $\beta_n^d$ follows directly from rearranging and simplifying the above formula for 
$\al_n^d$. Moreover, the resulting formula for $\phi(n)$ follows from the fact that if $d=1$, then $\al_n^1=n-1$ and 
$\beta_n^1=\phi(n)$. This concludes the proof of $(i)$.

For the proof of $(ii)$, note that if $n=p^t$, then the set of distinct primes of $n$ is $D_n=\{p\}$. Since the only nonempty subset of $\{p\}$ is itself and $\phi(p^t)=p^{t-1}(p-1)$, then $(ii)$ follows directly from $(i)$.

To prove $(iii)$, note that if  $n=p_1p_2$, then the only proper nonempty subsets of $D_n=\{p_1,p_2\}$ are $\{p_1\}$,
 $\{p_2\}$, and $\{p_1,p_2\}$. Moreover, $n_{p_1}=n/p_1=p_2$, $n_{p_2}=n/p_2=p_1$, and $n_{\{p_1,p_2\}}=n/p_1p_2=1$. Thus, it follows from $(i)$ that 
\begin{align*}
\al_{p_1p_2}^d
&=\beta_{p_1p_2}^d+ (-1)^{2} \al_{p_1}^d+ (-1)^{2} \al_{p_2}^d+ (-1)^{3} \al_{1}^d=\beta_{p_1p_2}^d+\al_{p_1}^d+\al_{p_2}^d,
\end{align*}
where we used $\al_{1}^d=0$ by Fact~\ref{fact:1}. 
\end{proof}

\subsection{Determining $\al_n^d$ for $d>n/2$ or $d\leq3$} \label{action}

We first determine $\al_n^d$ for $d>n/2$. We will use the following important theorem due to Savchev and Chen.

\begin{theorem}[\cite{SC}]\label{thm:SC}
Every zero-sum free $d$-tuple $\x$ in $\intsn^{d}$ of length $d>n/2$ can be uniquely represented as
 $(x_{1}k, x_{2}k, \dots, x_{d}k)$, where $k$ generates $\intsn$ and $x_{1}, x_{2}, \cdots, x_{d}$ are positive integers whose sum is less than $n$.
\end{theorem} 

We can now prove our main theorem in this section. 

\begin{theorem}\label{savchev-chen}
Let $n \ge 3$ and $d > n/2$. Then the number of zero-sum free $d$-tuples in $\ints_{n}$ is given by 
\[ \alpha_{n}^{d} = \phi(n) {n-1 \choose d}.\]
\end{theorem}

\begin{proof}
We consider the action of $Aut(\ints_{n})$ on $\G_n^{d}$. Since $d>n/2$, it follows from Theorem~\ref{thm:cha}$(d)$ that $\al_n^d=\phi(n)\cdot N$, where $N$ is the number of orbits under the action of $Aut(\ints_{n})$. So it suffices to determine $N$. Pick some zero-sum free $d$-tuple $\x$ in some orbit $O$. By Theorem~\ref{thm:SC},
$\x$ can be uniquely represented as $\x=(x_{1}k, x_{2}k, \dots, x_{d}k$), where $k$ generates $\intsn$ and $x_{1}, x_{2}, \cdots, x_{d}$ are positive integers such that $\sum_{i=1}^nx_i<n$. Note that if $\ell\in Aut(\ints_{n})$, then $\ell\cdot \x =\ell k(x_{1}, x_{2}, \dots, x_{d})$, where $\ell k\in Aut(\ints_{n})$ is another generator of $\intsn$. In fact $\ell k$ can be assumed to be any generator 
 of $\intsn$ by a suitable choice of $\ell\in Aut(\ints_{n})$. Thus, $N$ is the number of ordered tuples $(x_1,\ldots,x_d)$
 that satisfy $\sum_{i=1}^nx_i<n$, which is equivalent to the number of ordered partitions of $j$, with $1\leq j\leq n-1$, 
 into $d$ positive integers:
 \[N=\sum_{j=1}^{n-1}{j-1\choose d-1} ={n-1 \choose d}, \]
 where the last equality can be easily shown (e.g., by induction). Therefore, we obtain 
\[\al_{n}^{d} =\phi(n)\cdot N=\phi(n) {n-1 \choose d}.\]
\end{proof}

Note that a very similar theorem was proven by \cite{VP} for $d>\frac{2n}{3}$, where the author takes into account only minimal zero-sum free sequences as opposed to minimal zero-sum free tuples. 

\begin{corollary}
For any fixed positive integer $k$, we have $\alpha_n^{n-k} = \phi(n) {n-1 \choose k-1}$ for all large enough value of $n$. Moreover, 
\[ \liminf_n \frac{\alpha_{n+1}^{n+1-k}}{\alpha_n^{n-k}} = 0 , \ \ \text{ and } \ \ \limsup_n \frac{\alpha_{n+1}^{n+1-k}}{\alpha_n^{n-k}} = \infty.\]
\end{corollary}

\begin{proof}
For all large $n$, $n-k > n/2$, so by the above theorem we have
\[\alpha_n^{n-k} = \phi(n) {n-1 \choose n-k} = \phi(n) {n-1 \choose k-1}.\]
Therefore,
\[ \frac{\alpha_{n+1}^{n+1-k}}{\alpha_n^{n-k}} = \frac{\phi(n+1) {n \choose k-1}} {\phi(n) {n-1 \choose k-1}} = \frac{\phi(n+1)}{\phi(n)} \frac{n}{n-k+1}.\]

Since $n/(n-k+1)$ goes to $1$ as $n$ goes to $\infty$, the corollary now follows from a result of Somayajulu~\cite{SB} which states that 
\[ \underset{n}{\text{lim inf}} \ \ \frac{\phi(n+1)}{\phi(n)} = 0, \ \ \text{ and } \ \ \underset{n}{\text{lim sup}} \ \ \frac{\phi(n+1)}{\phi(n)} = \infty.\] 
\end{proof}

It is not hard to get explicit formulas for $\alpha_n^d$ for small values of $d$. For instance, $\alpha_n^1= n-1$ and 
$\alpha_n^2 = (n-1)(n-2)$. In principle, one can compute $\alpha_n^d$ for any $n$ and $d$ using the exclusion-inclusion
principle. We illustrate this when $d=3$ in the next proposition.
\begin{proposition}
Let $n \ge 3$ be a positive integer. Then 
\[ \alpha_n^3 = n^3 - 7n^2 + 15n-10+ \frac{1+(-1)^{n-1}}{2}\text{ for all } n \ge 3.\]
\end{proposition}
\begin{proof}  
Using the exclusion-inclusion principle, we compute the cardinality of the complement of $\G_n^3$ inside $(\ints_n \setminus \{ 0\})^3$. For any set $P$ of hyperplanes, we let $V(P)$ denote the subset of all points in $(\ints_n \setminus \{ 0\})^3$ which belong to the intersection of all planes in $P$. Then we have 
\[ \alpha_n^3 = (n-1)^3 - |V(x+y = 0) \cup V(y+z = 0) \cup V(x+z= 0) \cup V(x+y+z = 0)|.\] 
Now we compute cardinalities of various intersections. Note that $|V(x+y = 0)| = (n-1)^2$ because $x \ne 0$ and $z \ne 0$ are the only restrictions here, and $y = -x$. Similarly, $|V(y+z = 0)| = |V(x+ z = 0)| = (n-1)^2$. But $|V(x+y+z = 0) | = (n-1)(n-2)$ because while $x \ne 0$, $y \ne 0$ or $-x$ to ensure $z \ne 0$. We then have a total of $3(n-1)^2 + (n-1)(n-2)$ points.
Similarly, the reader can check that we have a total of $3(n-1)$ points coming from $2$-fold intersections.
 
 Now we look at $3$-fold intersections. As before, if any one of these intersections contains the plane $x+y+z=0$, then 
$x=0$, or $y=0$, or $z=0$. So we ignore these intersections. There is then only one intersection to consider, namely $V(x+y = 0, y+z = 0, z+x = 0) $. Adding all these equations gives $2(x+y + z)=0$. Since we are working in $\ints_n$, when $n$ is odd, this equation is equivalent to $x+y + z = 0$, so we can ignore this case. But when $n$ is even, $2(x+y + z = 0)$ would imply $x+y+z = n/2$. This in conjunction with $x+y = 0$ yields $z = n/2$. Similarly using $y+z = 0$ we get $x = n/2$, and using $x + z = 0$ we get $y = n/2$. Thus, when $n$ is even, we pick an additional point $(n/2, n/2, n/2)$. The total number of points in 3-fold intersections is $\theta_n$, which is $1$ if $n$ is even and $0$ otherwise.
 
Finally, a 4-fold intersection will give $(0,0,0)$. Now packing all these cardinalities in the exclusion-inclusion we get
\[ \alpha_n^d = (n-1)^3 - \left(3(n-1)^2 + (n-1) (n-2) - 3(n-1) +\theta_{n} - 0\right). \]
Simplifying this expression gives the desired formula.
\end{proof}

It will be very tedious to take this approach to compute $\alpha_n^d$ for higher values of $d$. However, using this approach, in the next section, we will show the existence of a polynomial of degree $d$ which captures $\alpha_n^d$ for suitable values of $n$.


\section{Hyperplane arrangements} \label{hyperplane}
A real hyperplane arrangement $\mathcal{H}$ is a finite collection of hyperplanes (subspaces of dimension $n-1$) in $\mathbb{R}^n$. 
An important polynomial attached to any real hyperplane arrangement is its characteristic polynomial; see \cite{SR} for the definition. This polynomial contains valuable information about $\mathcal{H}$. For instance, a fundamental result of Zaslavsky~\cite{ZT} states that the value of the characteristic polynomial at $-1$ gives, up to sign, the number of regions in a hyperplane arrangement. Similarly, the value at $1$ gives, up to sign, the number of bounded regions.  

 Now we look at hyperplane arrangements in $\mathbb{F}_p^d$ ($d$-dimensional vector space over $\mathbb{F}_p$). Note that $\G_p^d$, the space of zero-sum-free $d$-tuples, represents points in the complement of the hyperplanes in $\mathbb{F}_p^d$ defined by the $2^d - 1$ linear equations $\sum_{i \in S} x_i = 0$, where $S$ ranges over all the non-empty subsets of $\{ 1, 2, \cdots, d \}$.The following result is well-known, but we give independent proof and connect to the Hadamard bound. 
\begin{theorem} \label{prime-theorem}
Let $d$ be a positive integer. Then there exists a monic polynomial $f_d(x)$ of degree $d$ with integer coefficients such that for any sufficiently large prime $p$, $f_d(p) = \alpha_p^d$.
\end{theorem}
\begin{proof}
Consider the $2^d-1$ hyperplanes defined by the equations $\sum_{i \in S} x_i = 0$, where $S$ ranges over all nonempty subsets 
$S \subseteq \{ 1, 2, \cdots, d\}$. For every nonempty subset $P_i$ of these planes, let $A_i$ be a matrix whose rows are the coefficients of the linear equations that define the hyperplanes in $P_i$. Note that $A_i$ defines a linear map 
\[ A_i \colon\;\mathbb{F}_p^d \rar \mathbb{F}_p^{m_i},\]
where $m_i$ is the number of hyperplanes in $P_i$, and that $\text{Null}(A_i)$ is the intersection of hyperplanes on $P_i$. Then we have
\[ \alpha_p^d =p^d - |\text{Set of all points on the $2^d - 1$ hyperplanes}|. \]
We compute the cardinality of the union using exclusion-inclusion principle. We get an expression of the form
 \[ \alpha_p^d =p^d - \sum \pm |N(A_i)| .\]
By the rank-nullity theorem, $\text{Null}(A_i)$ is a $(d- \rank{A_i})$ dimensional subspace of $\mathbb{F}_p^d$. 
Therefore, $|\text{Null}(A_i)| =p^{d- \rank{A_i}}$.
This gives, 
 \[ \alpha_p^d =p^d - \sum \pm p^{d- \rank{A_i}}, \]
 where $ \rank{A_i} > 0 $ because $A_i$ are non-zero matrices. We will be done if we can show that $ \rank{A_i} $ does not depend on $p$ for a sufficiently large prime $p$. To this end, we use a result in linear algebra which states that the rank of a matrix $M$ is the maximal order of a non-zero minor of $M$. Since each $A_i$ is a matrix with 1s and 0s of a fixed order $m_i \times d$, there is a finite set of minors of all $A_i$'s over $\ints$. Let $v_d$ be the maximum prime which divides one of these minors. As soon as the $q$ exceeds $v_d$, $ \rank{A_i} $ becomes independent of $p$. 
 Then we have
 \begin{equation}\label{vd}
 \alpha_p^d =p^d - \sum \pm p^{d- \rank{A_i}}, \ \ p > v_d.
 \end{equation} 
 The RHS is clearly a polynomial in $p$ of degree $d$.
\end{proof}

By the general theory of hyperplane arrangements (see~\cite{SR}), it is known that for sufficiently large primes $p$, 
$\alpha_p^d = h_d(p)$ where $h_d(x)$ is the characteristic polynomial of the corresponding real hyperplane arrangement. 
Since the polynomial $f_d(x)$ which we constructed and $h_d(x)$ are both of degree $d$ and they agree for all sufficiently 
large primes, they are the same. These characteristic polynomials are considered to be intractable, in general. Here are the 
first few polynomials which are obtained using SageMath:
$f_1(x) = x-1$, $f_2(x) = x^2-3x+2 $, $f_3(x) = x^3-7x^2+15x-9$, $f_4(x) = x^4 -15x^3 + 80x^2-170x + 104$.

Having shown the existence of polynomials, which captures $\alpha_p^d$ for sufficiently large primes, two natural questions arise. 

\begin{enumerate}
\item What are the coefficients of these polynomials?
\item What is the value of the number $v_d$ shown in~\eqref{vd}? Can we get some bounds for it?
\end{enumerate}

 Recall from Section~\ref{sec:introduction} that $H_d$ denotes the $d \times (2^d-1)$ matrix whose columns are all the non-zero binary vectors in $\mathbb{F}_p^d$. 
Thus, $H_d$ is the matrix corresponding to our hyperplane arrangement. For instance, 
\[ H_3 = 
\left(\begin{array}{rrrrrrr}
1 & 0 & 0 & 1 & 1 & 1 & 1 \\
0 & 1 & 0 & 1 & 0 & 0 & 1 \\
0 & 0 & 1 & 0 & 1 & 1 & 1  
\end{array}\right).
\]
From the proof of the above theorem we see that the coefficient of $x^i$ ($0 \le i \le d-1$) of our polynomial $f_d(x)$ is obtained from all possible null spaces of dimension $i$ in the alternating sum given by the exclusion-inclusion principle. This gives:

\begin{proposition}\label{prop:coeff} The coefficient of $x^i$ ($0 \le i \le d-1$) of our polynomial $f_d(x)$ is given by $\sum_{j=1}^{2^d-1} \;(-1)^j m(j, i)$, where $m(j, i)$ is the number of subsets of $j$ columns of $H_d$ that span a $d-i$ dimensional subspace of
 $\mathbb{F}_p^d$.
\end{proposition}

A closed formula for these coefficients seems hard to obtain with the exception of one case. The coefficient of $x^{d-1} $ is given by $\sum_{j=1}^{2^d-1} \;(-1)^{j} m(j, d-1)$. Note that $m(1, d-1) = 2^d -1$ (every column is non-zero and there are $2^d-1$ columns) and $m(i, d-1) = 0$ for all $i > 1$ (a subset of size more than 1 cannot span a 1-dimensional subspace when working with binary vectors). Thus the coefficient of $x^{d-1}$ is $-(2^d -1)$. 

The constant term is $\sum_{j=1}^{2^d-1} \;(-1)^j m(j, 0)$. Note that $m(j,0)$ is the number of subsets of columns of $H_d$ of size $j$, which span $\mathbb{F}_p^d$. 

Now we consider the problem of finding a bound for $v_p$. 
\begin{proposition} \label{charprimes} Let $d$ be a positive integer. Then,
$\alpha_p^d = f_d(p) \; \text{ for all } \; p > d^{d/2}$.
In other words, $v_d\leq d^{d/2}$.
\end{proposition}
\begin{proof}
 Recall that if our prime $p$ exceeds the value of all possible minors obtained from $H_d$, then $\alpha_p^d = H_d(p)$. Now note that when computing all possible minors of $H_d$, it is enough to look at $d \times d$ submatrices because the minors from $j \times j$ $(j \le d)$ submatrices will be picked by those of $d \times d$ submatrices as can be seen by simply adding appropriate $0$s and $1$s. The maximum possible minor of a $d \times d$ submatrix of $H_d$ can be bounded by Hadamard's inequality which states that for any matrix $A$, $\det(A) \le \prod_{i=1}^d ||A_i||$, where $||A_i||$ is the Euclidean norm of the $ith $ column of $A$. Applying this inequality to a $d \times d$ submatrix of $M$ of $H_d$ gives 
 \[ \det(M) \le \prod_i ||M_i|| \le \prod_i \sqrt{d} = d^{d/2}.\] 
 This means that for all $p > d^{d/2}$ we have $\alpha_p^d = f_d(p)$.
 \end{proof}
 
 The maximum possible minor of a $d \times d$ submatrix of $H_d$ will be at most the maximum possible determinant of an arbitrary $0-1$ matrix of order $d$. This latter is called the Hadamard maximal determinant problem.
An upper bound for this determinant is $2^{-d}(d+1)^{(d+1)/2}$ (e.g., see~\cite{JB}). This shows that 
\[ \alpha_p^d = f_d(p) \mbox{ for all $p$ such that} \; p > \frac{(d+1)^{(d+1)/2}}{2^d}.\]
See entry A003432 in OEIS for the maximum possible determinant of a $(0,1)$-matrix of order $d$.

The following theorem offers a generalization of Theorem~\ref{prime-theorem}.

\begin{theorem}\label{thm:prime-ext}
For every positive integer $d$ there exists a monic polynomial $f_d(x)$ of degree $d$ with integer coefficients such that $\alpha_n^d = f_d(n)$ for all $n$ that are relatively prime to the determinant of any $d \times d$ binary matrix.
\end{theorem}
\begin{proof}
We start exactly as in the proof of Theorem \ref{prime-theorem}, but we work over $\ints_n$. With the same notation as there, we obtain
\[ \alpha_n^d = n^d - |\text{Set of all points on the $2^d - 1$ hyperplanes}|.\]
We compute the cardinality of the union using the exclusion-inclusion principle, and that gives
 \[ \alpha_n^d = n^d - \sum \pm |Ker(A_i)| .\]
 Note that since we are no longer working over a prime field, we look for the cardinality of the kernel of the group homomorphism 
 \[A_i \colon\; \mathbb{Z}_n^d \rar \mathbb{Z}_n^{m_i}, m_i \le d.\]
 Since $n$ is relatively prime to the determinant of any $d \times d$ binary matrix, in particular, it will also be relatively prime to any minor of $A_i$. These minors are exactly the collection of all the scalars that we multiply with in the process of converting $A_i$'s into their RRE forms over $\mathbb{Z}_n$. This will show that the under the given condition on $n$, $|Ker(A_i)| = n^{d-\rank(A_i)} $, where rank is computed over the reals. 
Thus, we have a polynomial function $f_d(x) = x^d - \sum \pm x^{d- \rank(A_i)}$, such that $\alpha_n^d = f_d(n)$ for all $n$ that are relatively prime to the determinant of any $d \times d$ binary matrix.
\end{proof}
\begin{corollary} \label{gcdlemma}
For every positive integer $d$, we have $\alpha_n^d = f_d(n)$ if $\text{gcd}\left(n, \lceil d^{d/2}\rceil!\right) = 1$. In particular, $\alpha_p^d = f_d(p)$ for all sufficiently large primes $p$.
\end{corollary}
\begin{proof}
Since the determinant of any binary $d \times d$ matrix is at most $d^{d/2}$, if we choose $n$ such that $\text{gcd}\left(n, \lceil d^{d/2}\rceil!\right)= 1$, then $n$ will be relatively prime to the determinant of any $d \times d$ binary matrix and we can apply Theorem~\ref{thm:prime-ext}.
\end{proof}
\section{Hypotheses} \label{conjectures}\label{sec:conj}

We examine the rows and columns of the following table of the values of $\alpha_n^d$ and state some hypotheses. The data is obtained using a program (available at~\cite{Sa}) that we wrote using SageMath software. The entries in the empty cells of Table~\ref{tab1} are zeros since $\al_n^d=0$ if $d\geq n$ (see Fact~\ref{fact:1}).
\begin{table}[ht]\label{tab1}
\begin{center}
\resizebox{\textwidth}{!}{
\begin{tabular}{|c|c|c|c|c|c|c|c|c|c|c|c|c|c|c|c|c|c|}
\hline
\diagbox{$n$}{$d$} &1&2&3&4&5&6&7&8&9&10&11&12&13&14&15&16&17 \\ \hline 	
2&1& & & & & & & & & & & & & & & & \\ \hline  	
3&2& 2 & & & & & & & & & & & & & & &\\ \hline 	
4&3& 6& 2 & & & & & & & & & & & & & &\\ \hline 
5&4&12 & 16& 4& & & & & & & & & & & & &\\ \hline 	
6&5& 20& 44 & 10 & 2 & & & & & & & & & & & &\\ \hline 
7&6& 30& 96& 90& 36& 6& & & & & & & & & & &\\ \hline 
8&7& 42& 174& 240& 84& 28& 4& & & & & & & & & &\\ \hline 
9&8& 56& 288& 690& 336&168 & 48& 6 & & & & & & & & &\\ \hline 
10&9&72&440&1344 & 984& 336& 144& 36& 4 & & & & & & & &\\ \hline 
11&10& 90& 640& 2590& 3060& 2100& 1200& 450& 100& 10 & & & & & & &\\ \hline 
12&11& 110& 890& 4330& 5786& 2436& 1320& 660& 220& 44& 4 & & & & & & \\ \hline 
13&12& 132& 1200& 7020& 14832& 12264& 9504& 5940& 2640& 792& 144& 12 & & & & & \\ \hline 
14&13& 156& 1572& 10560& 26172& 22686& 13992& 7722& 4290& 1716& 468& 78& 6 & & & &\\ \hline 
15&14& 182& 2016& 15564& 52488& 49392& 28736& 24024& 16016& 8008& 2912& 728& 112& 8 & & &\\ \hline 
16&15& 210& 2534& 21840& 83292& 95620& 73876& 56880& 40040& 24024&10920& 3640& 840& 120& 8 & &\\ \hline 
17&16& 240& 3136& 30160& 143616& \mbox{\bf 217056}& \mbox{\bf 208000}& \mbox{\bf 209808}& {\bf 183040}& \mbox{\bf 128128}& 69888& 29120& 8960& 1920& 256& 16 &\\ \hline 
18&17& 272& 3824& 40330& 217574& {\bf 326088}& 292080& 216672& \mbox{\bf 162780}& \mbox{\bf 116688}& 74256& 37128& 14280& 4080& 816& 102& 6\\ \hline 
\end{tabular}
}
\end{center}
\caption{$\alpha_n^d$ for $2\leq n\leq 18$ and $1\leq d\leq n-1$.}
\label{tab2}
\end{table}  	

 Analyzing the data in Table~\ref{tab2} led us to consider the following hypotheses.
 \vskip 3mm

\noindent
 \textbf{Row hypothesis:} For a fixed $n$, $\{ \alpha_n^d \}$ is an increasing sequence in $d$ for $d < n/2$ and a decreasing sequence in $d$ for $d > n/2$.

Note that for $d > n/2$, $\alpha_n^d = \phi(n) {n-1 \choose d}$. Since, for a fixed value of $n$, the binomial coefficients ${n-1 \choose d}$ form a decreasing sequence for $d > n/2$, this proves the second part of the hypothesis. The first part of the hypothesis is false. The smallest counterexample is when $n = 17$: $\alpha_{17}^6 > \alpha_{17}^7 $. We expect that there would be infinitely many counter-examples.

 \vskip 3mm
\noindent
\textbf{Column hypothesis:} For a fixed $d$, $\{ \alpha_n^d \}$ is an increasing sequence in $n$.

A quick check in SageMath reveals many counterexamples to the inequality involving the Euler-phi function. 
For instance, we have a counterexample at $(n, d) = (17, 6)$, since $\alpha_{18}^6 < \alpha_{17}^6$.
To investigate this hypothesis, we will use the following lemma. 
\begin{lemma}\label{lem:col}
For $d > (n+1)/2$, we have $\alpha_{n+1}^d > \alpha_{n}^d$  if and only if $\frac{\phi(n+1)}{\phi(n)} > \frac{n-d}{n}$.
\end{lemma}
\begin{proof}
Recall that $\alpha_n^d = \phi(n) {n-1 \choose d}$ for $d > n/2$. From this, we have 
$\alpha_{n+1}^d >\alpha_{n}^d$ is equivalent to $\phi(n+1) {n \choose d} > \phi(n) {n-1 \choose d}$, which 
in turn is equivalent to $\frac{\phi(n+1)}{\phi(n)} > \frac{n-d}{n}$.
\end{proof}
We now show that the column hypothesis fails for infinitely many columns.
\begin{proposition}\label{prop:nocol}
There are infinitely many positive integers $d$ such that $\alpha_{n+1}^d < \alpha_n^d$ for some $n$ that depends on $d$. 
\end{proposition}
\begin{proof}
A result of Somayajulu~\cite{SB} states that 
\[ \underset{n\to\infty }{\text{lim inf}} \ \ \frac{\phi(n+1)}{\phi(n)} = 0. \] 
In particular, there are infinitely many positive integers $n_i$ such that 

\[\ \frac{\phi(n_i+1)}{\phi(n_i)} \le \frac{1}{4}. \]
Since $ \frac{1}{4} = 1 - \frac{3}{4} = 1 - (\frac{3n_i/4}{n_i} ) \le 1- \left(\frac{[ 3n_i/4 ]}{n_i}\right)$, 
we obtain 
\[\ \frac{\phi(n_i+1)}{\phi(n_i)} \le \frac{n_i - [ 3 n_i/4 ]}{n_i}.\]

Setting $d_i =[ 3 n_i/4 ] $, we see from the above lemma that the last inequality is equivalent to 
\[ \alpha_{n_i + 1}^{d_i} < \alpha_{n_i}^{d_i},\]
because $d_i = [ 3 n_i/4 ] > (n_i+1)/2$. 
This shows that the column hypothesis fails in infinitely many columns. 
\end{proof}
 \vskip 3mm

The failure of the column hypothesis leads to a revised hypothesis.
\vskip 3mm
\noindent
\textbf{Eventual column hypothesis:} For any fixed value of $d$, the sequence $\alpha_n^d$ is eventually an increasing sequence in $n$, i.e., $\alpha_n^d$ is an increasing sequence if $n$ is large enough.
\vskip 3mm
Note that exact formulas for $\alpha_n^d$ for $d \le 3$, and the existence of a monic polynomial $f_d(x)$ of degree $d$, such that $\alpha_n^d = f_d(n)$, whenever $\text{gcd}\left(n, \lceil d^{d/2}\rceil!\right) = 1$, give support to this hypothesis.  In fact, in Section \ref{sec:asymptotic}, we will prove the eventual column hypothesis. 

\section{Bounds}\label{sec:bnd}
In this section, we obtain some upper bounds and lower bounds for $\alpha_n^d$ and $\beta_n^d$ and use them in the next section in conjunction with the characteristic polynomials and number-theoretic results to study the asymptotic behavior of $\alpha_n^d$ as a sequence in $n$ when $d$ is fixed. It is enough to restrict to interesting cases: $n \ge 3$ and $d \ge 3$.

It is easy to obtain upper bounds for $\alpha_n^d$ using simple counting techniques, as illustrated in the subsequent two propositions.  

\begin{proposition} We have 
 $ 0 \le \beta_n^d \le \alpha_n^d \le (n-1)^{d-1}(n-2)$ for $d \ge 3$ and all $n \ge 3$. 
\end{proposition}

\begin{proof}
The relation $0 \le \beta_n^d \le \alpha_n^d$ is clear by definition. For the upper bound, let $(x_1, \cdots, x_d)$ be an element in $\G_n^d$. Then we know that for all $1 \le i \le d-1$, we have $x_i \ne 0$, and $x_d \ne 0$ or $-(x_1 + \cdots +x_{d-1})$. So we have at most $n-2$ choices for $x_d$ and $n-1$ choices for the rest. 
\end{proof}

One can obtain a slightly sharper upper bound by restricting to primes.
 
\begin{proposition}
For all primes $p \ge 3$ and all $d \ge 3$, $\beta_p^d = \alpha_p^d \le (p-1) (p-2)^{d-2}(p-3)$. 
\end{proposition}

\begin{proof}
In $\ints_p^d$, every zero-sum-free $d$-tuple is also irreducible because the condition $\gcd(x_1, \cdots, x_d, p) =1$ is automatically satisfied. 
For $1 \le i \le p-1$, let $\G_p^d(i)$ be the set of all $d$-tuples in $\G_p^d$ where the first component is $i$. Since $p$ is a prime, multiplication by $i$ induces a bijection between $\G_p^d(1)$ and $\G_p^d(i)$.
Thus, we have $|\G_p^d| = (p-1) |\G_p^{d}(1)|$. Note that 
\[ \G_p^{d}(1) \subseteq \{ (1, x_2, \cdots, x_d) : x_i \ne 0 \text{ or } -1\}.\]
We have at most $p-2$ choices for each $x_i$. In addition, $x_d$ cannot be equal to $- (1 + x_2 + \cdots+ x_{d-1})$, and this quantity cannot be $0$ or $-1$ when $p \ge 3$. This completes the proof.
\end{proof}

We now turn to lower bounds. The following observation can be used to get some recursive lower bounds for $\alpha_n^d$. Let $m$ be a divisor of $n$. Then the natural ring homomorphism $\mathbb{Z}_n \rightarrow \mathbb{Z}_m$ extends to a ring homomorphism $\ {\psi} \colon\, \mathbb{Z}_n^d \rightarrow \mathbb{Z}_m^d$.
It is clear that $\mathbf{x}$ is in $\G_n^d$ whenever ${\psi} (\mathbf{x})$ is in $\G_m^d$.

\begin{proposition} For all $m$ and $n$ such that $m$ divides $n$, we have 
$\alpha_n^d \ge \left(\frac{n}{m}\right)^d \alpha_m^d $
\end{proposition}

\begin{proof} The kernel of the homomorphism $\psi \colon\, \mathbb{Z}_n^d \rightarrow \mathbb{Z}_m^d $
has order $(n/m)^d$. (In fact, if $(m)$ is the ideal generated by $m$ in $\mathbb{Z}_n$, then the kernel is $(\mathbb{Z}_n/(m))^d$). So every zero-sum free $d$-tuple in $\mathbb{Z}_m^d$ pulls back to $(n/m)^d$ zero-sum free $d$-tuples in $\mathbb{Z}_n^d$. This gives the stated lower bound.
 \end{proof}

To get better and explicit lower bounds, we use results from Section~\ref{action}. Recall that $\I_n^d$ is the set of all $d$-tuples $\x\in\G_n^d$ such that $\gcd(x_1,\ldots,x_n,n)=1$.

\begin{proposition} \label{bounds} For any integers $n$ and $d$ such that $1\leq d \le n-1$, we have 

\n $(i)$ $\al_n^d \geq \binom{n-1}{d}$ and, in particular, $\alpha_n^d \ge \frac{(n-d)^d}{d!}$;

\n $(ii)$ $\beta_n^d \ge \phi(n)\binom{n-2}{d-1}$.  

\end{proposition}
\begin{proof}
Part $(i)$ follows from the observation that every $d$-tuple $(x_1,\ldots,x_d)$ of positive integers
 that satisfies $\sum_{i=1}^dx_i<n$ is clearly zero-sum-free. The cardinality of the latter was shown in the proof of Theorem \ref{savchev-chen} to be ${n-1 \choose d}$. Thus,
\[\al_n^d\geq \binom{n-1}{d}=\frac{(n-1)\ldots (n-d)}{d!} \geq \frac{(n-d)^d}{d!}.\]

To prove part $(ii)$, let $\x=(1,x_2,\ldots,x_d)$ be a $d$-tuple of positive integers such that $x_2+\ldots+x_d\leq n-2$.
Then $\x\in\G_n^d$ and, therefore, $k\x\in\G_n^d$ for any $k$ such that $\gcd(k,n)=1$. Since $\gcd(k,n)=1$ implies that 
$\gcd(k,kx_2,\ldots,kx_d,n)=1$, 
it follows that $k\x\in\I_n^d$. The number of $(d-1)$-tuple $(x_2,\ldots,x_d)$ such that $x_i\geq 1$ and 
$x_2+\ldots+x_d\leq n-2$ is equal to the number of ordered partitions of $j$, with $1\leq j\leq n-2$, into $d-1$ 
positive integers is:
\[\sum_{j=1}^{n-2}{j-1\choose d-2} ={n-2 \choose d-1}. \]
Thus, 
\[\beta_n^d \ge |\{\x=(k,kx_2,\ldots,kx_d):\, \gcd(k,n)=1,\,\, 0<x_2+\ldots+x_d\leq n-2\}|=\phi(n)\binom{n-2}{d-1}.\]
\end{proof}

\section{Asymptotic Results} \label{sec:asymptotic}

In this section we will focus on asymptotic results on the sequences $\{ \alpha_n^d \}$ and $\{\beta_n^d \}$. We say that $a_n = \mathcal{O} (b_n)$ if there is a positive constant $K$ and an integer $N$ such that $|a_n| \le K b_n$ for all $n \ge N$. 
 We begin with the sequence $\{ \alpha_n^d \}$.

\begin{theorem}\label{thm:asymptotic-m}
Let $d$ be a fixed positive integer. 

\n $(i)$ We have $\alpha_n^d = n^d - (2^d-1)n^{d-1} + \mathcal{O}(n^{d-2})$.

\n $(ii)$  The sequence $\{\alpha_n^d\}$ is asymptotically equivalent to the sequence $\{n^d\}$. That is, 
$\lim\limits_{{n}\to\infty}\; \frac{\alpha_n^d}{n^d} =1$.

\n $(iii)$ We also have $\lim\limits_{{n}\to\infty}\frac{\alpha_{n+1}^d}{\alpha_n^d} = 1$.
\end{theorem}
\begin{proof}
We work within the framework of the proof of Theorem \ref{prime-theorem}, but over $\mathbb{Z}_n$, instead of $\mathbb{F}_p$. Consider the $2^d-1$ hyperplanes defined over $\mathbb{Z}_n$ by the equations $\sum_{i \in S} x_i = 0$, where $S$ ranges over the non-empty subsets of $\{1, 2, \cdots, d\}$. Let $\mathcal{P}$ be the collection of all non-empty subsets of these hyperplanes. Corresponding to each collection $P$ of hyperplanes, we consider a matrix $A_P$ whose rows are the coefficients of the linear equations which define the hyperplanes in $P$. We can view $A_p$ as a group homomorphism,
\[A_P \colon\, \mathbb{Z}_n^d \rightarrow \mathbb{Z}_n^{m_P}, \]
where $m_P = | P|$. Then,
\[ \bigcap_{H \in P} H = \ker(A_P).\]
Applying the exclusion-inclusion to compute the number of points in the union of the $2^d-1$ hyperplanes gives the following: 
\begin{eqnarray*}
\alpha_n^d & = & n^d - |\text{union of the } 2^d-1 \text{ hyperplanes}|\\
& = & n^d - \left( \sum_{P \in \mathcal{P}} (-1)^{|P|+1} | \bigcap_{H \in P} H | \right) \\
& = & n^d - \left( \sum_{P \in \mathcal{P}} (-1)^{|P|+1} |\ker(A_P)| \right)\\
& = & n^d - \left( \sum_{P \in \mathcal{P},\, |P|=1 } |\ker(A_P) | + \sum_{P \in \mathcal{P},\, |P|\ge 2} (-1)^{|P|+1}|\ker(A_P)| \right).
\end{eqnarray*} 
We now claim that
\[ |\ker(A_P)| =
\begin{cases}
n^{d-1} & \text{ if } |P| = 1,\\
\le n^{d-2} & \text{ if } |P| \ge 2.\\
\end{cases}
\]
When $|P| = 1$, $|\ker{A_P}| = |\ker([a_1, \cdots, a_d])|$, where $[a_1, \cdots, a_d]$ are coefficients of the hyperplane in question. Since all coefficients are $0$ or $1$, and not all $0$, it is clear that $|\ker{A_P}| = n^{d-1}$. Similarly when $|P| \ge 2$, note that $\ker{A_P}$ is contained in the kernel of a $2 \times d$ matrix obtained by taking the first two rows of $A_P$. The cardinality of the latter is $n^{d-2}$ because the two rows are distinct binary vectors.

Substituting these values of $|\ker(A_P)|$ in the above expression for $\alpha_n^d$ and using the fact that there are $2^d-1$ hyperplanes in our collection and $|\mathcal{P}| = 2^{2^d-1} -1$, we get  
\begin{equation}\label{eq:asymp}
\alpha_n^d = n^d - (2^d-1)n^{d-1} + \mathcal{O}(n^{d-2}),
\end{equation}
which proves $(i)$.

Finally, it follows from~\eqref{eq:asymp} that
\[
\lim_{n\to \infty} \frac{\alpha_n^d}{n^d}=\lim_{n\to \infty} \frac{n^d - (2^d-1)n^{d-1} + \mathcal{O}(n^{d-2})}{n^d}=1,\]
and
\[\lim_{n\to \infty} \frac{\alpha_{n+1}^d}{\alpha_n^d} =\lim_{n\to \infty} \frac{(n+1)^d - (2^d-1)(n+1)^{d-1} + \mathcal{O}(n^{d-2})}{n^d - (2^d-1)n^{d-1} + \mathcal{O}(n^{d-2})}=1,\]
proving part $(i)$ and part $(ii)$, respectively.
\end{proof}

We also have the following application of Theorem \ref{thm:asymptotic-m}.

\begin{theorem}\label{thm:eventual-col}
The eventual column hypothesis is true. That is, for any fixed positive integer $d$, the sequence $\{ \alpha_n^d\}$ is an increasing sequence if $n$ is large enough.
\end{theorem}
\begin{proof}
Using the formula for $\alpha_n^d$ from part $(i)$ of Theorem \ref{thm:asymptotic-m} and the binomial theorem, we have
\begin{eqnarray*}
&&\alpha_{n+1}^d - \alpha_n^d \\
& = & \left[(n+1)^d - (2^d-1)(n+1)^{d-1} + \mathcal{O}(n^{d-2}) \right] - 
\left[ n^d - (2^d-1)n^{d-1} + \mathcal{O}(n^{d-2}) \right] \\
 & = & [ (n+1)^d - n^d]-(2^d -1) \left[ (n+1)^{d-1} - n^{d-1} \right] + \mathcal{O}(n^{d-2})\\
 & = & \left[ dn^{d-1} + \mathcal{O}(n^{d-2}) \right] - (2^d-1) \mathcal{O}(n^{d-2}) + \mathcal{O}(n^{d-2})\\
 & = & dn^{d-1} + \mathcal{O}(n^{d-2}).
 \end{eqnarray*}
In particular, this shows that for all sufficiently large values of $n$, $\alpha_{n+1}^d - \alpha_n^d > 0$. This completes the proof of the theorem. 
\end{proof}

We now turn our attention to $\beta_n^d := |\I_n^d| $. What can be said about the growth of the sequence? We do not expect this sequence to be asymptotic to $n^d$. In fact, if $d = 1$, then $\beta_n^1 = \phi(n)$, where it is known that 
\[ 0 = \liminf_n \frac{\phi(n)}{n} < \limsup_n \frac{\phi(n)}{n} = 1. \]

To analyze this sequence, we use the following theorem of Hardy and Wright~\cite[Page 267, Theorem 327]{HW}. For any $\epsilon > 0$, 
\[\lim\limits_{n\to\infty} \frac{\phi(n)}{n^{1- \epsilon}} = \infty.\]
In other words, the order of $\phi(n)$ is nearly equal to $n$ when $n$ is large enough.

\begin{theorem}\label{thm:limsup}
Let $d$ be a fixed positive integer. Then, 
\begin{enumerate}
\item $\{ \beta_n^d/n^d \}$ is a sequence in $[0, 1]$ with  $\limsup_n \frac{\beta_n^d}{n^d} \;= 1$;
\item $\beta_n^d$ is $\mathcal{O} (n^d)$ but not $\mathcal{O} (n^{d- \epsilon})$ for any $\epsilon > 0$. 
\end{enumerate}
 \end{theorem}
 \begin{proof}
Since $0 \le \beta_n^d \le n^d$, it is clear that the sequence in question belongs to $[0, 1]$.
To compute $\limsup$, we use the fact that whenever $\{ x_n\}$ is a sequence in $[a, b]$ with a subsequence $\{ x_{n_k}\}$ whose limit is $b$, then $\limsup_n x_n = b$. So it is enough to show that there is a subsequence of $\{ \frac{\beta_n^d}{n^d} \}$ whose limit is 1. To this end, we consider the subsequence $\{ \frac{\beta_p^d}{p^d} \}$ that corresponds to primes and note that $\beta_p^d = \alpha_p^d$, and for a fixed $d$, there is a monic polynomial $f_d(x) = x^d + a_{d-1}x^{d-1} + \cdots + a_1 x + a_0$ of degree $d$ such that $\alpha_p^d = f_d(p)$ for all sufficiently large primes. Then we have
 \[ \limsup_n \frac{\beta_n^d}{n^d} =\lim\limits_{p \to\infty} = \frac{\beta_p^d}{p^d} = \lim\limits_{p \to\infty} \frac{f_d(p)}{p^d} = \lim\limits_{p \to\infty} \frac{p^d + a_{d-1}p^{d-1} + \cdots + a_1 p + a_0}{p^d} = 1.\]
 
The fact that $\beta_n^d$ is $\mathcal{O} (n^d)$ follows from the trivial upper bound $\beta_n^d \le n^d$ because this bound is a polynomial in $n$ of degree $d$. Let $\epsilon > 0$ be fixed. Suppose to the contrary that 
$\beta_n^d$ is $\mathcal{O} (n^{d- \epsilon})$. Then, by definition, there is a constant $K$ and an integer $N$ such that 
\[ \beta_n^d \le K n^{d - \epsilon} \mbox{ for all }  n \ge N.\]
Using the lower bound for $\beta_n^d$ from Proposition \ref{bounds} (i), we get 
\[ \phi(n)\binom{n-2}{d-1} \le \beta_n^d \le K n^{d - \epsilon} \mbox{ for all }  n \ge N. \]
Simplifying this gives
\[ \frac{\phi(n)}{n^{1-\epsilon}} \le K \frac{n^{d-1}}{ \binom{n-2}{d-1} } \mbox{ for all }  n \ge N. \]
Using the above-mentioned result of Hardy and Wright, the LHS goes to infinity as $n$ goes to infinity. However, the RHS tends to $K(d-1)!$ as $n$ goes to infinity. This contradiction shows that $\beta_n^d$ is not $\mathcal{O} (n^{d- \epsilon})$ for any $\epsilon > 0$.
\end{proof}

What about $\liminf \beta_n^d$? Although the answer for $d=1$, as noted above, is 0, for $d \ge 2$, we will see that the answer is greater than 0. To show this, we first recall some results from analytic number theory.  

The Riemann zeta function is an important complex-valued function in number theory that is an analytic continuation of the series defined by 
\[ \zeta(s) = \prod_{p \text{ prime } }\left(1 - \frac{1}{p^s}\right)^{-1} = \sum_{n=1}^{\infty} \frac{1}{n^s}, \mbox{ for all }  \text{Re}(s) > 1.\]
This function is the key to understanding the distribution of primes and many number theoretic problems. 
For our problem, we will use the following result.  

The probability that $d$ $(\ge 2)$ randomly chosen positive integers will have a gcd of 1 is given by \cite{NJ, TA}
\[\ \prod_{p \text{ prime }} \left(1- \frac{1}{p^d} \right) = \frac{1}{\zeta(d)}.\]
Note that a point $(v_1, v_2, \cdots, v_d)$ in $\mathbb{Z}^d \setminus \{ \bar{0} \}$ is visible from the origin if and only if $\gcd(v_1, \cdots, v_d) = 1$.
In particular, when $d = 2$, this probability is $1/\zeta(2) = 6/\pi^2$ -- the density of lattice points on the plane that are visible from the origin.

\begin{proposition}\label{prop:zeta}
For any positive integer $d \ge 2$, $\liminf_n \frac{\beta_n^d}{n^d} \ge \frac{1}{d! \zeta(d)}.$
\end{proposition}
\begin{proof}
Recall that $\beta_n^d$ is the cardinality of all irreducible zero-sum-free $d$-tuples in $\mathbb{Z}_n^d$. That is, $\beta_n^d = | \mathcal{I}_n^d | $, where 
\[ \mathcal{I}_n^d = \mathcal{G}_n^d \cap \{ (x_1, \cdots, x_d) \in \mathbb{Z}_n^d \colon \gcd(x_1, \cdots, x_d, n) = 1 \}.\]
Note that 
\[ \{ (x_1, \cdots, x_d) \in (\mathbb{Z}_n \setminus \{ 0 \})^d \colon\, x_1 + \cdots+ x_d < n \} \subseteq \mathcal{G}_n^d, \]
and moreover, $\gcd(x_1, \cdots, x_d) = 1$  implies that $\gcd(x_1, \cdots, x_d, n) = 1$. 
Since we identify elements of $\ints_n$ with the representatives $\{ 0, 1, \cdots, n\}$ of nonnegative integers, we have,
\[ \{ (x_1, \cdots, x_d) \in \mathbb{Z}_{>0}^d \; \colon \; 0 < x_1 + \cdots+ x_d < n \text{ and } \gcd(x_1, \cdots, x_d) = 1 \} \subseteq \mathcal{I}_n^d.\]
Let $R_n^d$ denote the region in $\mathbb{R}_{> 0}^d$ that is enclosed by the coordinate planes and the hyperplane $x_1 + x_2 + \cdots +x_d = n$. 
Using multivariable calculus it can be shown that
\[ \text{Vol}(R_n^d) = \int_0^n \int_0^{n- x_1} \int_0^{n- x_1 -x_2 } \cdots \int_0^{n - x_1 -x_2 -\cdots -x_{d-1}} dx_d\; dx_{d-1} \cdots dx_1 = \frac{n^d}{d!}. \]
Let $\theta_n^d$ be the number of lattice points in the interior of the region $R_n^d$ that are visible from the origin.
The above inclusion shows that 
\[
 \frac{\theta_n^d}{\text{Vol}(R_n^d) } \frac{1}{d!}  
 = \frac{\theta_n^d}{\text{Vol}(R_n^d) } \frac{\text{Vol}(R_n^d) }{n^d} = \frac{\theta_n^d}{n^d} \le \frac{\beta_n^d}{n^d} . 
\]
As $n$ goes to infinity, $R_n^d$ goes to $\mathbb{R}_{> 0}^d$, and therefore $\frac{\theta_n^d}{\text{Vol}(R_n^d) }$ goes to $1/\zeta(d)$. So taking $\liminf$ on both sides of the above inequality gives
\[\frac{1}{\zeta(d)} \frac{1}{d!} \le \liminf_n \frac{\beta_n^d}{n^d}. \]
 This completes the proof.
\end{proof}

Since $\zeta(d) > 0$ for $d\ge 2$, the following theorem is now clear from the above lower bound and the definition of $\liminf$. This gives another proof of part (B) of Theorem \ref{thm:limsup}.
\begin{theorem}\label{thm:bounded}
Let $d \ge 2$ be a positive integer. Then there is a positive integer $N_d$ such that 
\[ \left(\frac{1}{d!\; \zeta(d)} \right) n^d \; \le \; \beta_n^d \; \le \; n^d \mbox{ for all }  n > N_d. \]
In other words, the sequence $\{ \beta_n^d \}$ is asymptotically bounded above and below by the sequence $\{ n^d \}$. 
\end{theorem}

\section{Appendix: Mathieu-Zhao subspaces} \label{Mathieu}

In this final section, we explain how we arrived at computing the number of zero-sum-free sequences. 
Let $A$ be a commutative $k$-algebra, for a field $k$. A $k$-subspace $M$ of $A$ is said to be a Mathieu-Zhao subspace if the following property holds: Let $a$ belong to $A$ be such that $a^m \in M$ for all $m \ge 1$. Then for any $b \in A$, we have $ba^m \in M$ for all $m$ sufficiently large, i.e., there exists a positive integer $N \ge 1$ (that depends on both $a$ and $b$) such that $ba^m \in M$ for all $m \ge N$. This definition resembles the definition of an ideal in $A$. It is easy to check that every ideal in $A$ is a Mathieu-Zhao subspace, but the converse is not true. So a Mathieu-Zhao subspace can be viewed as a generalization of the concept of an ideal. They were introduced by Zhao \cite{ZW} in his study of the Jacobian conjecture. They have played a central role in the subsequent work related to the Jacobian conjecture and other related conjectures. The following idempotent criterion of Zhao gives a useful characterization of these subspaces. 

\begin{proposition}[\cite{ZW}]
Let $M$ be a $k$-subspace of a unital commutative finite dimensional $k$-algbera $A$. $M$ is a Mathieu-Zhao subspace of $A$ if and only if for every idempotent $e$ in $M$, the ideal $(e)$ is also contained in $M$.
\end{proposition}

We learned the following example of a Mathieu-Zhao subspace from Wenhua Zhao.
\begin{proposition}[\cite{ZW2}]\label{prop:mat}  
Let $p$ be an odd prime. The kernel of the linear functional defined by the vector $(c_1, c_{2}, \cdots, c_{n})$ in $\ints_{p}^{n}$ is a Mathieu-Zhao subspace
if and only if for any subset $S \subseteq \{1, 2, \dots, n\}$ either $\sum_{i \in S} c_{i} \ne 0$ or $c_{i} = 0$ for all $i$ in $S$.
\end{proposition}
\begin{proof}
This follows directly from the idempotent criterion. We begin by noting that any idempotent $\e$ in $\ints_p^{n}$ satisfies $e_i^2 = e_i$, which means $e_i = 0$ or $1$ for all $i$. Now let $\text{Ker} (\cv)$ denote the kernel of the linear functional defined by the vector $\cv= (c_1, c_2, \cdots, c_n)$. The kernel of $\cv$ is a Mathieu-Zhao subspace if and only if for any idempotent
 $\e$ in $\text{Ker} (\cv)$, $\x \e$ (component-wise multiplication) also belongs to $\text{Ker} (\cv)$ for all 
$\x\in \ints_p^n$. This means, for any $S \subseteq \{ 1, 2, \cdots, n\}$, $\sum_{i \in S}c_i = 0$ 
 implies that $\sum_{i \in S} x_i c_i = 0$ for all $x_i$. If some $c_i \ne 0$, this condition fails. So it must be the case that for any $S \subseteq \{ 1, 2, \cdots,n\}$, either $c_i = 0$ for all $i$ or $\sum_{i \in S} c_i \ne 0$
\end{proof}

The above theorem provides us with some examples of Mathieu-Zhao subspaces. We were interested in counting the number of Mathieu-Zhao subspaces from the above proposition. To this end, we say that a vector $(c_1, c_{2}, \cdots, c_n)$ in $\ints_{p}^n$ is Mathieu-Zhao if it satisfies the condition given in Proposition~\ref{prop:mat}. Because of that proposition, it is enough to count the number of Mathieu-Zhao $n$-tuples.

\begin{proposition} If $\mathcal{M}_n$ denote the collection of all Mathieu-Zhao $n$-tuples in $\mathbb{F}_{p}^n$, then
\[ |\mathcal{M}_{n}| = {n \choose 1} \alpha_{n}^1 + {n \choose 2} \alpha_{n}^2 + \cdots + {n\choose n} \alpha_{n}^{n}.\]
Thus, to compute $|\mathcal{M}_{n}|$, we need to compute $\alpha^d_{n}$ for ${1\le d \le n}$.
\end{proposition}

\begin{proof}
Let $T_{d}$ be the set of Mathieu-Zhao $n$-tuples with exactly $d$ non-zero entries. Then 
$\mathcal{M}_n$ is the disjoint union of subsets $T_d$, where for $1\leq d\leq n$ there are ${n \choose d}$ ways to select the positions for the $d$ non-zero entries. From the definitions, it is clear that $T_{d}$ is equal to ${n\choose d} \al_n^{d}$.
\end{proof}

\noindent {\bf Acknowledgements.}
We thank the referee for valuable comments which helped us improve the presentation of the paper. We also thank the editors for several editorial corrections that made the paper more readable.

We  thank Gail Yamskulna and Wenhua Zhao for the discussions which led to the problem of computing the number of zero-sum-free tuples. We also thank Richard Stanley for bringing to our attention the characteristic polynomials of hyperplane arrangements and their connection to our work. Some of this research was done when the first author was visiting MIT. He would like to thank MIT for its hospitality. The first author also presented this work at the University of Tokyo at a Hyperplane Arrangements and Singularities conference. He got valuable feedback from Christin Bibby, Graham Denham, and Alex Suiu. He thanks them for their input. 


\end{document}